\documentclass[a4paper,american,12pt]{amsart}
 
\usepackage[latin1]{inputenc}

\usepackage{amsmath, bbm}
\usepackage{amssymb}
\usepackage{hyperref}
\usepackage{url}
\usepackage{babel}
\usepackage{tikz}
\usepackage{pgf,tikz,pgfplots}
\usepackage{graphicx}
\usepackage{tkz-graph}

\usetikzlibrary{arrows}

\usetikzlibrary{calc}
\usepackage{amsfonts}
\input xy
\xyoption{all}
\newcommand{\cl}{\mathrm{cl}}
\newcommand{\trop}{\mathrm{trop}}
\newcommand{\supp}{\mathrm{supp}}
\newcommand{\rk}{\mathrm{rk}}
\newcommand{\Aut}{\mathrm{Aut}}
\newcommand{\crem}{\mathrm{Crem}}
\newcommand{\R}{\mathbb{R}}
\newcommand{\Z}{\mathbb{Z}}

\newtheorem{thm}{Theorem}[section]
\newtheorem{defi}[thm]{Definition}

\newtheorem{prop}[thm]{Proposition}

\newtheorem{lemma}[thm]{Lemma}

\newtheorem{corollary}[thm]{Corollary}

          {\theoremstyle{definition}
}
          {\theoremstyle{definition}
\newtheorem{example}[thm]{Example}
}

\newtheorem{ques}[thm]{Question}

\newtheorem*{corintro}{Corollary}

\title{Combinatorial Cremona automorphisms and Coxeter arrangement matroids }
\author{Stefan Rettenmayr and Annette Werner}

\begin{document}

\begin{abstract}
We explore birational geometry of matroids by investigating automorphisms of their coarse Bergman fans. Combinatorial Cremona maps provide such automorphisms of Bergman fans which are not induced by matroid automorphisms. We investigate the structure of matroids allowing combinatorial Cremona maps and prove a realizability criterion in the presence of two different Cremonas. We also prove that for all matroids associated to Coxeter arrangements the group of coarse automorphisms of the Bergman fan is generated by the matroid automorphisms and at most one combinatorial Cremona map. 
\end{abstract}

\maketitle
\markright{Combinatorial Cremona automorphisms and Coxeter arrangement matroids}
\centerline{MSC(2020): 14T20, 52B40, 14E07}

\section{Introduction}

In this paper we investigate automorphisms of Bergman fans of matroids endowed with the coarse fan structure. Some, but in general not all, such automorphisms arise from matroid automorphisms. For matroids which are realizable by an essential hyperplane arrangement in projective space, birational morphisms which are regular on the arrangement complement give rise to such automorphisms induced by invertible linear maps of the ambient space. Therefore we regard isomorphisms
of coarse Bergman fans as matroid analogs of birational maps. 

The rich geometry of realizable matroids tends to inspire purely combinatorial 
questions for general matroids. It is often a challenge to study which geometrical properties persist to all matroids.  
A famous recent example  is the development of Hodge theory of matroids by Adiprasito, Huh and Katz \cite{ahk}.

In \cite{sw}, combinatorial analogs of Cremona transformations were introduced. In some realizable examples, such maps are in fact tropicalizations of algebraic geometric Cremona transformations, see section \ref{sec:A} for the example of type $A$ root systems. 

Combinatorial Cremona maps on the Bergman fan of a simple, connected matroid $M$ have the property that there is a basis of $M$ such that its indicator vectors are mapped to indicator vectors of corank one flats (see Definition \ref{def:cremona}). 
The present paper is motivated by the open question:
\begin{ques}\label{ques} For which simple connected matroids which are not parallel connections in a non-trivial way is the group of automorphisms of the coarse Bergman fan generated by Cremonas and maps induced by matroid automorphisms?     
\end{ques}

After explaining some background notions on Bergman fans and their automorphisms in section \ref{sec2},
we discuss the example of the braid arrangement in section \ref{sec:A}. This illustrates the connection of combinatorial Cremona maps with classical Cremonas for realizable matroids. Moreover, it was shown in \cite{sw} that Question \ref{ques} has a positive answer for the braid arrangement matroid.

Section \ref{sec:rank3} deals with rank 3 matroids where question \ref{ques} has a positive answer for all simple connected matroids which are not parallel connections. 
In section \ref{sec:cremona} we look at matroids admitting a combinatorial Cremona map and investigate the structural consequences for the flats in $M$. 
We introduce the support graph of subsets of $M$ with respect to a Cremona basis as a new combinatorial tool.
Section \ref{sec:2cremonas} deals with matroids admitting at least two different combinatorial Cremona maps. Theorem \ref{twocremonas} shows that the support graph of one Cremona basis with respect to the other Cremona basis consists of star-shaped trees. We deduce that the matroid automorphisms act transitively on the set of Cremona bases:
\begin{corintro}[\bf \ref{cor: involution}] 
Assume that the matroid $M$ admits two Cremona bases $b$ and $b^\ast$. Then there exists an involutive matroid automorphism $\varphi\colon E(M) \rightarrow E(M)$ mapping $b$ to $b^\ast$.
\end{corintro}

We also  prove the following  realizability criterion in the presence of two different Cremona bases, see Theorem  \ref{thm:repcrit} for a more detailed  statement. 
\begin{thm}
Let $M$ be simple and connected of rank $d+1 \geq 3$ and assume that $M$ admits two different Cremona bases $b = \{b_0, \dotsc, b_d\}$ and $b^\ast = \{ b_0, b_1^\ast, \dotsc, b_d^\ast\}$ with $|b \cap b^\ast|=1$.
Then  $M$ is realizable over all fields with at least $|E(M)| - d$ elements. 
\end{thm}

In section \ref{sec:Cox}
 we determine the automorphism group of the coarse Bergman fan for all matroids associated to Coxeter arrangements. It turns out that only for arrangements of type $A_n$ and $B_n$
this group is strictly bigger than the group of matroid automorphisms. The $B_n$ case follows from its relation to Dowling matroids which were investigated in \cite{sw}. Combining Theorem \ref{aut:A} (for $A_n$), Proposition \ref{prop:B} (for $B_n$) and Theorem \ref{auto} we find that for all matroids associated to Coxeter arrangements the automorphism group of the coarse Bergman fan is generated by matroid automorphisms and combinatorial Cremona maps, which gives a positive answer to question \ref{ques}.


{\bf Notation:} In this paper, graphs are multigraphs, i.e.~we allow parallel edges. For a graph $G$, we denote by $V(G)$ the set of vertices, and by $E(G)$ the set of edges.

{\bf Acknowledgements: } We thank the referee for his or her useful comments. Research on this paper was partially supported by the Deutsche Forschungsgemeinschaft (DFG, German Research Foundation) TRR 326 \textit{Geometry and Arithmetic of Uniformized Structures}, project number 444845124.

\section{Matroids, Bergman fans and combinatorial Cremonas} \label{sec2}

A matroid is a pair $M = (E, r)$ where $E= E(M)$ is a finite set  and $r\colon \mathcal{P}(E) \to \Z_{\geq 0}$ is a rank function satisfying the axioms

\begin{enumerate}
\item $r(A) \geq 0 $ for all $A \subseteq E$ and $r(\emptyset) = 0$,
\item For $A, B \subseteq E$ we have $r(A \cup B) + r(A \cap B) \leq r(A) + r(B)$, 
\item For any $i \in E$ and $A \subseteq E$, we have $r(A) \leq r(A \cup \{ i \}) \leq r(A) + 1$. 
\end{enumerate}

The rank of a matroid $M$ is $r(E)$ and is denoted $r(M)$. 
A flat of a matroid $M$ is a subset $F \subseteq E$ which is closed under the rank function, which means that for all $i \notin F$ the rank of $F \cup \{ i\}$ is strictly greater than the rank of $F$. 

For every subset $S$ of $E$ the closure $\cl(S)$ of $S$ is the inclusion minimal flat containing $S$. If $F_1$ and $F_2$ are flats of $M$, the join $F_1 \vee F_2$ is defined as the closure of $F_1 \cup F_2$. 

Since we are interested in Bergman fans of matroids, we will generally assume that $M$ is simple. We denote by $(v_i)_{i \in E}$ the canonical basis of the real vector space $\mathbb{R}^E$ generated by the ground set $E$ of $M$. The (affine) Bergman fan $\tilde{B}(M)$ associated to the loopfree matroid $M$ is defined as the subset of $\R^{E}$ consisting of all vectors $\sum_i a_i v_i \in \mathbb{R}^E$ such that for each circuit $C$ of $M$ the minimum of the set $\{a_i: i \in C\}$ is attained at least twice.
It is invariant under scaling by the vector $\mathbbm{1} = \sum_i v_i$. We will mostly work with the \emph{projective Bergman fan} $B(M)$, which is defined as the image of the affine Bergman fan under the quotient map $\R^{E} / \R \mathbbm{1}$.
Note that for matroids realizable by a hyperplane arrangement $\mathcal{A}$ in projective space, the Bergman fan coincides with the tropicalization of the hyperplane complement $\Omega_{\mathcal{A}}$ with respect to the torus embedding given by the hyperplane coordinates.

The set $B(M)$ carries several fan structures. There is the fine fan structure, where the cones are given by all flags of flats in $M$. A coarser structure is the minimal nested set structure where rays are given by the connected flats of $M$, and cones correspond to so-called {nested} collections of connected flats. Here a collection $\mathcal{S}$  of connected flats is called nested if for any subcollection of incomparable flats $F_1, \dots, F_t \in \mathcal{S}$ with $t \geq 2$, the join is not connected. We will also use the coarse fan structure which is defined as the coarsest fan structure available on $B(M)$. If the criterion \cite[Theorem 5.3]{feistu} is satisfied, it agrees with the minimal nested set structure.

\begin{defi}{\cite[Definition 2.5]{sw}} \label{def:auto} Let $M$ and $M'$ be simple matroids  on ground sets $E$ and $E'$, respectively.
 An isomorphism of Bergman fans  $ B_\ast(M)$ and  $B_\ast(M')$, where $\ast$ denotes  choices of fan structures, is  a linear map $\phi\colon \R^E / \R \mathbbm{1}  \to \R^{E'}/ \R \mathbbm{1} $ derived from a lattice isomorphism  $\Z^E / \Z \mathbbm{1} \to \Z^{E'} / \Z \mathbbm{1}$, such that $\phi$ restricts to an isomorphism of fans  $ B_\ast(M) \to B_\ast(M')$.
 \end{defi}

Note that any matroid isomorphism between $M$ and $M'$ induces an isomorphism of their Bergman fans. 
For simple matroids $M$ and  $M'$ of rank $2$, every isomorphism between their Bergman fans (where all fan structures coincide) is obviously induced by a matroid isomorphism. Hence our matroids will generally have rank at least $3$. The rank $3$ case is discussed in the next section.

The study of automorphisms of the Bergman fan of a matroid with respect to different fan structures was initiated in \cite{sw}. It turns out that every isomorphism between the fine structures of Bergman fans is induced by a matroid automorphism, see \cite[Theorem 6.3]{sw}. However, this is no longer the case if one considers coarser fan structures. In \cite{sw} candidates for automorphisms of Bergman fans which are not induced by matroid automorphisms are introduced, inspired by birational geometry applied to realizable matroids. As explained in \cite[Section 3]{sw}, in the case of a matroid given by an essential hyperplane arrangement $\mathcal{A}$ in projective space, any dominant endomorphism of the hyperplane complement $\Omega_{\mathcal{A}}$ extends to an endomorphism of the intrinsic torus and thus tropicalizes to a unimodular linear map preserving the Bergman fan.
Since unimodularity holds in this case for realizable matroids, the previous definition for all matroids also stipulates unimodularity.

For general matroids, \cite[Definition 8.1]{sw}  gives the following purely combinatorial notion of an analog of a Cremona transformation for a simple matroid $M$ on the ground set $E$ with a basis $b = \{b_0, \dots, b_{d}\}$: 
Define the flat $B_j$ as the closure of $b \backslash \{ b_j\}$. For any subset $T \subseteq E$ we write 
$v_T = \sum_{k \in T} v_k \in \R^E$.

\begin{defi}{\cite[Definition 8.1]{sw}}
\label{def:cremona} We define the $\Z$-linear map $\crem_{b} \colon \R^{E} \to \R^{E}$ by 
$v_{b_j} \mapsto v_{B_j}$ for the basis elements $b_j$, and 
$v_k \mapsto v_k$ for all $k \in E \backslash b$. 
If $\crem_b$ maps the line $\R \mathbbm{1}$ to itself, we also write 
$\crem_b\colon \R^{E } / \R \mathbbm{1} {\rightarrow} \R^{E} / \R \mathbbm{1}$ for the quotient map, and call $\crem_b$ a combinatorial Cremona map. 
\end{defi}

In particular, a combinatorial Cremona map is defined as a unimodular linear map on the ambient space of the Bergman fan. It is shown in \cite[Theorem 8.3]{sw} that the map $\crem_b$ descends to a linear map $\crem_b\colon \R^{E} / \R \mathbbm{1} \rightarrow \R^{E} / \R \mathbbm{1}$ mapping $B(M)$ to itself, if only if the sets $\cl\{b_i, b_j\} \backslash \{b_i, b_j\}$ partition the set $E \backslash b$ into pairwise disjoint subsets.
We call any basis of a matroid $M$ with this property a \emph{Cremona basis}. Moreover, \cite[Lemma 8.4]{sw} shows that $\crem_b$ acts as an involution on $B(M)$.

For simple matroids $M$, it is an interesting challenge to investigate the structure of the automorphism group $\Aut(B_c(M))$ of the Bergman fan endowed with its coarse fan structure. 

Our motivating question \ref{ques} is to determine all simple connected matroids $M$ which are not parallel connections in a non-trivial way such that the automorphism group $\Aut(B_c(M))$ is generated by matroid automorphisms and combinatorial Cremona maps. 
If $M$ does not admit any Cremona basis, this boils down to the question whether or not $\Aut(M) \rightarrow \Aut(B_c(M))$ is surjective.

Note that if $M$ is a parallel connection, then its coarse Bergman fan is a product of two coarse Bergman fans (see e.g.~\cite[Theorem 7.2]{sw}).
In this case, a Cremona automorphism on one factor may induce an automorphism of the whole Bergman fan which is no longer a Cremona automorphism. Therefore question \ref{ques}  only deals with matroids which are not parallel connections in a non-trivial way.

By \cite[Theorem 1.2]{sw},  all simple modularly complemented matroids belong to the class of matroids discussed in question \ref{ques}.
In the present paper we show in section \ref{sec:Cox} that this is also true for all matroids given by a Coxeter arrangement.

\section{The braid arrangement}\label{sec:A}
As a motivating example we look at the matroid $M$ associated to the braid arrangement, i.e.~the hyperplane arrangement given by the root system of type $A_n$. 

Hence let $K$ be an algebraically closed field and consider the set of hyperplanes given by the root system of type $A_n$ for $n \geq 3$ in $K^{n+1}$. After passing to the quotient by their intersection, we obtain the essential arrangement $\mathcal{A}$ in $\mathbb{P}^{n-1}_{K}$ given by the hyperplanes $V(x_i)$ for $i = 0, \dotsc, n-1$ and $V(x_i - x_j)$ for $i < j$ in $\{0, \dotsc, n-1\}$.

The realizable matroid $M$ induced by this arrangement $\mathcal{A}$ is isomorphic to the graphic matroid given by the complete graph $K_{n+1}$. This implies that the automorphism group of $M$ is equal to the symmetric group $S_{n+1}$, i.e.~the Weyl group of the $A_n$ root system.  

It is well-known that the associated hyperplane complement $\Omega_{\mathcal{A}}$ can be identified with the moduli space $M_{0,n+2}$ of $n+2$-pointed curves of genus zero. 
The minimal wonderful compactification of $\Omega_{\mathcal{A}}$ is isomorphic to the Deligne-Mumford compactification $\overline{M}_{0, n+2}$. It is shown in \cite{brunomela} that the automorphism group of $\overline{M}_{0, n+2}$ is isomorphic to $S_{n+2}$ and hence equal to the automorphism group of $\Omega_{\mathcal{A}}$.

The Bergman fan of the matroid $M$ can be identified with the moduli space $M_{0,n+2}^{trop}$ of tropical genus zero curves with $n+2$ marked points, see \cite[Chapter 4]{ArdilaKlivans},  so it has automorphism group given by $S_{n+2}$ (see \cite{ap}). 

This difference between matroid automorphisms and automorphisms of the coarse Bergman fans can be explained by the presence of combinatorial Cremona automorphisms acting as new transpositions. In fact, we have the following result: 

\begin{thm}[{\cite[Proposition 3.2]{sw}}] \label{aut:A} The automorphism group of the coarse Bergman fan $B_c(M)$ is generated by the matroid automorphisms of $M$ and one combinatorial Cremona map. 
\end{thm}

Let us have a closer look at such combinatorial Cremona maps for $n=3$:

\begin{example}\label{ex:A3} The matroid associated to the Coxeter arrangement of type $A_3$ has 6 elements. In the presentation introduced above, we can realize it by the essential hyperplane arrangement $\mathcal{A}$ in $\mathbb{P}_K^2$ given by $x=0$, $y=0$, $y=z$, $x=z$, $x=y$ and $z=0$. We can identify it with the matroid associated to the complete graph on four vertices with edges labeled by $1, \dotsc, 6$ (see below) by using the order of this list of elements. Then $M$ has four Cremona bases, which correspond to the star-shaped spanning trees in the following figure: 

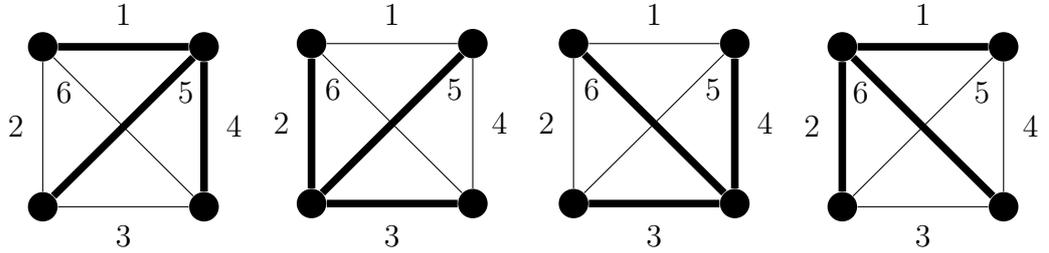
\begin{center}
    \begin{tikzpicture}
        \node[style={circle, fill}] (1) at (45:1.5){};
        \node[style={circle, fill}] (2) at (135:1.5){};
        \node[style={circle, fill}] (3) at (225:1.5){};
        \node[style={circle, fill}] (4) at (315:1.5){};
        \draw[line width=1mm] (1) -- (2) node[midway, above=3pt] {1};
        \draw (2) -- (3) node[midway, left=3pt] {2};
        \draw (3) -- (4) node[midway, below=3pt] {3};
        \draw[line width=1mm] (4) -- (1) node[midway, right=3pt] {4};
        \draw[line width=1mm] (1) -- (3) node[near start, right=2pt] {5};
        \draw (2) -- (4) node[near start, left=2pt] {6};
    \end{tikzpicture}
    \begin{tikzpicture}
        \node[style={circle, fill}] (1) at (45:1.5){};
        \node[style={circle, fill}] (2) at (135:1.5){};
        \node[style={circle, fill}] (3) at (225:1.5){};
        \node[style={circle, fill}] (4) at (315:1.5){};
        \draw (1) -- (2) node[midway, above=3pt] {1};
        \draw[line width=1mm] (2) -- (3) node[midway, left=3pt] {2};
        \draw[line width=1mm] (3) -- (4) node[midway, below=3pt] {3};
        \draw (4) -- (1) node[midway, right=3pt] {4};
        \draw[line width=1mm] (1) -- (3) node[near start, right=2pt] {5};
        \draw (2) -- (4) node[near start, left=2pt] {6};
    \end{tikzpicture}
    \begin{tikzpicture}
        \node[style={circle, fill}] (1) at (45:1.5){};
        \node[style={circle, fill}] (2) at (135:1.5){};
        \node[style={circle, fill}] (3) at (225:1.5){};
        \node[style={circle, fill}] (4) at (315:1.5){};
        \draw (1) -- (2) node[midway, above=3pt] {1};
        \draw (2) -- (3) node[midway, left=3pt] {2};
        \draw[line width=1mm] (3) -- (4) node[midway, below=3pt] {3};
        \draw[line width=1mm] (4) -- (1) node[midway, right=3pt] {4};
        \draw (1) -- (3) node[near start, right=2pt] {5};
        \draw[line width=1mm] (2) -- (4) node[near start, left=2pt] {6};
    \end{tikzpicture}
    \begin{tikzpicture}
        \node[style={circle, fill}] (1) at (45:1.5){};
        \node[style={circle, fill}] (2) at (135:1.5){};
        \node[style={circle, fill}] (3) at (225:1.5){};
        \node[style={circle, fill}] (4) at (315:1.5){};
        \draw[line width=1mm] (1) -- (2) node[midway, above=3pt] {1};
        \draw[line width=1mm] (2) -- (3) node[midway, left=3pt] {2};
        \draw (3) -- (4) node[midway, below=3pt] {3};
        \draw (4) -- (1) node[midway, right=3pt] {4};
        \draw (1) -- (3) node[near start, right=2pt] {5};
        \draw[line width=1mm] (2) -- (4) node[near start, left=2pt] {6};
    \end{tikzpicture}
\end{center}

For example, the combinatorial Cremona automorphism $\crem_b$ associated to the Cremona basis $b = \{ 1, 2, 6 \}$ is given by the following linear automorphism of $\R^6$: 
\[ v_1 \mapsto v_{\{2, 3, 6\}}, ~ v_2 \mapsto v_{\{1, 4, 6\}}, ~ v_3 \mapsto v_3, ~ v_4 \mapsto v_4, ~ v_5 \mapsto v_5, ~ v_6 \mapsto v_{\{1, 2, 5\}}. \] 

The following figure illustrates the cone complex of the Bergman fan of $M$ endowed with the minimal nested set structure. In this figure, vertices correspond to  rays of the cone complex, and  edges are given by adjacency of flats. The Cremona automorphism $\crem_b$ acts by rotating the outer hexagon by 180 degrees.

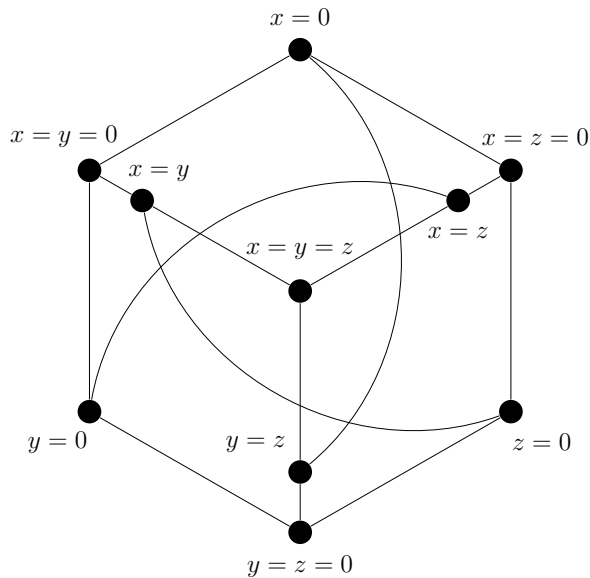
\begin{center}
\scalebox{0.8}{
\begin{tikzpicture}[every node/.style={circle,fill}]
	\node[label={[label distance=-10pt]90:{$x=0$}}] (1) at (90:4){};
	\node[label={[label distance=-17pt]150:{$x=y=0$}}] (2) at (150:4){};
	\node[label={[label distance=-5pt]210:{$y=0$}}] (3) at (210:4){};
	\node[label={[label distance=-20pt]270:{$y=z=0$}}] (4) at (270:4){};
	\node[label={[label distance=-5pt]330:{$z=0$}}] (5) at (330:4){};
	\node[label={[label distance=-17pt]30:{$x=z=0$}}] (6) at (30:4){};
	\draw (1) -- (2) -- (3) -- (4) -- (5) -- (6) -- (1);
	\node[label={[label distance=-12pt]10:{$x=y$}}] (10) at (150:3){};
	\node[label={[label distance=-30pt]150:{$y=z$ \quad \quad \phantom{a} }}] (11) at (270:3){};
	\node[label={[label distance=-10pt]270:{$x=z$}}] (12) at (30:3){};
	\node[label={[label distance=-17pt]90:{$x=y=z$}}] (13) at (0:0){};
    \draw[bend left=50] (1) to (11);
    \draw[bend left=50] (3) to (12);
    \draw[bend left=50] (5) to (10);
    \draw (2) -- (10) -- (13);
    \draw (4) -- (11) -- (13);
    \draw (6) -- (12) -- (13);
\end{tikzpicture}}
\end{center}

Note that this graph is isomorphic to the Petersen graph, see also \cite[Example 3.3]{feistu}.

The combinatorial Cremona automorphism $\crem_b$ is the tropicalization of the standard Cremona transformation $\crem\colon \mathbb{P}_K^2 \to \mathbb{P}_K^2$, $[x : y : z] \mapsto [yz : xz : xy]$, since the Cremona basis $b$ corresponds to the coordinate hyperplanes $\{ x = 0, y = 0, z = 0 \}$. More precisely, the coordinate map of the arrangement
\[ j\colon \mathbb{P}_K^2 \to \mathbb{P}_K^5, \quad [x : y : z] \mapsto [x : y : y-z : x-z : x-y : z] \]
identifies the hyperplane complement $\Omega_\mathcal{A} \subset \mathbb{P}_K^2$ with the linear subvariety $V(x_0 - x_1 - x_4, x_0 - x_5 - x_3, x_1 - x_5 - x_2)$ of the torus $T = \mathbb{G}_{m,K}^6 / \mathbb{G}_{m,K} \subset \mathbb{P}_K^5$ with diagonally embedded $\mathbb{G}_{m,K}$.  If we define the monomial map $\varphi\colon T \to T$ by
\[ [ x_0 : x_1 : x_2 : x_3 : x_4 : x_5 ] \mapsto [ x_1x_5 : x_0x_5 : -x_0x_2 : -x_1x_3 : -x_4x_5 : x_0x_1 ], \]
then the diagram
\[
\xymatrix{
\Omega_\mathcal{A} \ar[d]_\crem \ar@{^{(}->}[r]^j & T \ar[d]^\varphi & [ x : y : z ] \ar@{|->}[r] \ar@{|->}[d] & [x : y : y-z : x-z : x-y : z ] \ar@{|->}[d] \\
\Omega_\mathcal{A} \ar@{^{(}->}[r]^j & T & [ yz : xz : xy ] \ar@{|->}[r] & [ yz : xz : x(z-y) : y(z-x)  : z(y-x) : xy ]
}
\]
commutes. Thus the tropicalization $\mathrm{trop}(\varphi): \mathbb{R}^6 / \mathbb{R}\mathbbm{1} \rightarrow \mathbb{R}^6 / \mathbb{R}\mathbbm{1}$
 is the linear map given by multiplication with the matrix
\[ \begin{pmatrix} 0 & 1 & 0 & 0 & 0 & 1 \\ 1 & 0 & 0 & 0 & 0 & 1 \\ 1 & 0 & 1 & 0 & 0 & 0 \\ 0 & 1 & 0 & 1 & 0 & 0 \\ 0 & 0 & 0 & 0 & 1 & 1 \\ 1 & 1 & 0 & 0 & 0 & 0 \end{pmatrix} \]
and hence coincides with $\crem_b$.

\end{example}

\section{The rank 3 case}\label{sec:rank3}
For simple matroids of rank $3$ which are not parallel connections in a non-trivial way, \cite[Theorem 9.5]{sw} tells us that every  automorphism  of $B_c(M)$ is either induced by a matroid automorphism or the composition of a Cremona map  (in the sense of Definition \ref{def:cremona}) and a matroid automorphism. 

Note that by Definition \ref{def:auto}, every automorphism of $B_c(M)$ is obtained by restricting a linear automorphism of the ambient space which respects a lattice.

Less restrictively, one could also look at arbitrary linear automorphisms of the ambient space preserving the Bergman fan with, say, the minimal nested set fan structure. These automorphisms give rise to automorphisms of the associated minimal nested set complex of the matroid. In this way, we can interpret the natural duality for arrangements consisting of all rational hyperplanes in projective spaces over finite fields, as the following example of the Fano matroid shows. 

\begin{example}
    The Fano matroid is the simple rank 3 matroid $F_7$ on the ground set $E = \{ 1, \dotsc, 7\}$ with rank 2 flats $\{ 1, 2, 4 \}$, $\{ 1, 3, 6 \}$, $\{ 1, 5, 7 \}$, $\{ 2, 3, 5 \}$, $\{ 2, 6, 7 \}$, $\{ 3, 4, 7 \}$, and $\{ 4, 5, 6 \}$.

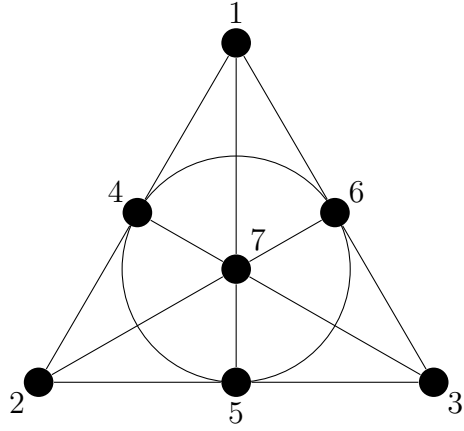
    \begin{center}
        \begin{tikzpicture}[every node/.style={circle,fill}]
        	\node[label={[label distance=-5pt]90:{$1$}}] (1) at (90:3){};
        	\node[label={[label distance=-5pt]150:{$4$}}] (2) at (150:1.5){};
        	\node[label={[label distance=-5pt]210:{$2$}}] (3) at (210:3){};
        	\node[label={[label distance=-5pt]270:{$5$}}] (4) at (270:1.5){};
        	\node[label={[label distance=-5pt]330:{$3$}}] (5) at (330:3){};
        	\node[label={[label distance=-5pt]30:{$6$}}] (6) at (30:1.5){};
        	\draw (1) -- (2) -- (3) -- (4) -- (5) -- (6) -- (1);
        	\node[label={[label distance=-2pt]80:{$7$}}] (7) at (0:0){};
            \draw (1) -- (7) -- (5);
            \draw (2) -- (7) -- (6);
            \draw (3) -- (7) -- (4);
            \draw (0:0) circle (1.5);
        \end{tikzpicture}
    \end{center}
    
    The self-duality of the Fano plane in the sense of projective geometry (not to be confused with duality of matroids) induces a linear automorphism $\phi\colon \R^E/\R \mathbbm{1} \to \R^E/\R \mathbbm{1}$ given by
    \begin{align*}
        &v_1 \mapsto v_{\{2, 3, 5\}}, ~ v_2 \mapsto v_{\{1, 3, 6\}}, ~ v_3 \mapsto v_{\{1, 2, 4\}}, ~ v_4 \mapsto v_{\{3, 4, 7\}} \\
        &v_5 \mapsto v_{\{1, 5, 7\}}, ~ v_6 \mapsto v_{\{2, 6, 7\}}, ~ v_7 \mapsto v_{\{4, 5, 6\}},
    \end{align*}
    which preserves the minimal nested set complex. However, the linear map $\phi$ is realized with respect to the basis $\{ v_1, \dotsc, v_6 \}$ by the matrix
    \[ \begin{pmatrix} 0 & 1 & 1 & -1 & 0 & -1 \\ 1 & 0 & 1 & -1 & -1 & 0 \\ 1 & 1 & 0 & 0 & -1 & -1 \\ 0 & 0 & 1 & 0 & -1 & -1 \\ 1 & 0 & 0 & -1 & 0 & -1 \\ 0 & 1 & 0 & -1 & -1 & 0 \end{pmatrix} \]
    whose determinant is $-8$, hence $\phi$ is not induced from a lattice isomorphism.
\end{example}


\section{Matroids admitting a Cremona basis}\label{sec:cremona}
Let $M$ be a simple and connected matroid of rank $d+1 \geq 3$ on the ground set $E = E(M)$. 
In this section we study matroids such that the associated coarse Bergman fan admits a Cremona automorphism in the sense of Definition \ref{def:cremona}.
Recall that by \cite[Theorem 8.3]{sw}, this is equivalent to the existence of a basis $b = \{b_0, \ldots, b_d\}$  of $M$ such that the 
sets
\[F_{ij}= \cl\{b_i,b_j\} \backslash \{b_i, b_j\}\]
partition the set $E \backslash b$ into pairwise disjoint subsets. Note that $F_{ij}$ may be empty. We call any basis $b$ of $M$ with this property a Cremona basis of $M$.

\begin{defi} Assume that $M$ has a Cremona basis $b$, and let $S$ be a subset of $E$.
We define the $b$-support of $S$ as the set \[\supp_b(S) = (b \cap S) \cup \bigcup_{F_{ij} \cap S \neq \emptyset} \{b_i, b_j\}. \]
\end{defi}
Note that $\rk(S) \leq |\supp_b(S)|$ since $S \subseteq \cl(\supp_b(S))$. For the next definition recall that our graphs are multigraphs. 


\begin{defi}
Assume that $M$ has a Cremona basis $b$. 
\begin{enumerate}
    \item 
For every subset $S$ of $E$ we define the support graph $G_b(S)$ of $S$ with respect to $b$ as the graph with vertex set $\supp_b(S)$ such that every element in the set $S \cap F_{ij}$ gives rise to an edge between $b_i$ and $b_j$. 
\item
We call $S$ support-connected if the graph $G_b(S)$ is connected.
\end{enumerate}
\end{defi}

\begin{lemma}\label{typeofflats} Assume that $M$ has a Cremona basis $b$, and consider a non-empty flat $F$ of $M$ with connected support graph $G_b(F)$.
Then one of the two following conditions holds:
\begin{enumerate}
\item $F \cap b = \supp_b(F)$, in which case we call $F$ a basis flat, or
\item $F \cap b = \emptyset$, in which case we call $F$ a non-basis flat.
\end{enumerate}
\end{lemma}
\begin{proof} 
Assume that $F \cap b$ is non-empty, and consider any element $b_i \in F \cap b$. Let $b_j \in b$ be a  vertex in $\supp_b(F)$ adjacent to the vertex $b_i$. Then there exists an edge $e \in F \cap F_{ij}$ between $b_i$ and $b_j$, which implies that $\{ b_i, b_j, e \}$ is a circuit as $M$ is simple. Therefore $b_j$ is also contained in $F \cap b$. Since $G_b(F)$ is connected, we conclude that $F \cap b$ is equal to $\supp_b(F)$.\end{proof}

Note that every connected flat is support-connected. Indeed, if $H$ is a connected component of $G_b(F)$, then $F$ is contained in the direct sum of $\cl(V(H))$ and $\cl(b \setminus V(H))$.
Note further that a basis flat $F$ satisfies $\rk(F) = |F \cap b|$, since it is generated by $\supp_b(F)$.

\begin{lemma}\label{nonbasisflats} Assume that $M$ has a Cremona basis $b$ as above, and 
let $F$ be a non-basis flat of $M$ (with connected support graph $G_b(F)$.)  

\begin{enumerate} 
\item For all $i,j \in \{0, \ldots, d\}$ the intersection $F \cap F_{ij} $ is either empty or consists of one element only. 
\item $F$ contains at most $\binom{|\supp_b(F)|} {2}$ elements. 
\item If $F$ is connected of rank $2$, the support graph $G_b(F)$ is a triangle.
\item We have $\rk(F)= |\supp_b(F)| -1$.

\end{enumerate}
\end{lemma}

\begin{proof}
If $F \cap F_{ij}$ contains two different elements $e,f$, then $\cl\{e,f\}$ is a flat of rank two contained in $\cl\{b_i,b_j\}$, since $M$ is simple. Hence both flats are equal. This implies that $F \cap b \neq \emptyset$, so that $F$ cannot be a non-basis flat. This shows the first claim. The second claim follows immediately from the first. In order to show the third claim, assume that $F$ has rank $2$, so that we can write $F = \cl\{e,f\}$. Since $e$ and $f$ do not lie in the same $F_{ij}$, but generate a connected flat, we find that the support graph of $F$ has three vertices. By our first claim, the support graph of $F$ is a simple graph. On the other hand, as $F$ is connected, the support graph contains at least one edge beyond $e$ and $f$, so that $G_b(F)$ must indeed be a triangle. In order to prove the fourth statement, let $b_i$ be an element in $\supp_b(F)$. Then $F$ is a strict subset of the flat $\cl(F \cup \{b_i\})$, which is a basis flat  with the same support as $F$. This implies our claim. \end{proof}

\section{Matroids with two Cremonas}\label{sec:2cremonas}
The following result shows the relation between two different Cremona bases. Note that the matroid given by the Coxeter arrangement associated to $A_n$ discussed in section \ref{sec:A} provides an example of a coarse Bergman fan admitting several Cremona automorphisms.

\begin{thm}\label{twocremonas}
Let $M$ be a simple and connected matroid of rank $d+1 \geq 3$  which has two different Cremona bases $b$ and $b^*$.
If $H$ is a connected component of $G_b(b^\ast)$, then $V(H) \cap b^\ast = \{ b_0 \}$ for some $b_0 \in b$, and the component $H$ is simple and star-shaped with center $b_0$, i.e. its $|V(H)|-1$ edges connect $b_0$ to all other elements in $V(H)$.
In particular, the support graph $G_b(b^\ast)$ is connected if and only if $|b \cap b^\ast| = 1$. 

Moreover, if $b_i, b_j  \in b$ are vertices in two different components in $G_b(b^\ast)$ such that $F_{ij} \neq \emptyset$, then $b_i$ and $b_j$ lie in $b \cap b^\ast$.
\end{thm}

\begin{proof} We begin with the following observation. Let $H$ be any connected component of the support graph $G_b(b^\ast)$. Put $I = \cl(V(H)) \cap b^\ast$. Then $H$ is the support graph of the (independent) subset $I$ of $ b^\ast$ with respect to the basis $b$. As $H$ is connected, we have $|E(H)| \geq |V(H)| -1$. Since $|V(H)| = |\supp_b(I)| \geq \rk(I) = |I|$ by Lemma \ref{typeofflats} and Lemma \ref{nonbasisflats}, this implies that $|I \cap b |\leq 1 $.

Consider first the case that $G_b(b^\ast)$ is connected, i.e. $H = G_b(b^\ast)$ and $I = b^\ast$.
Then we have just seen that $|b^\ast \cap b| \leq 1$.   Since $b^\ast$ is a Cremona basis, every element $b_i$ in $b\backslash b^\ast$ is contained in some $\cl\{b^\ast_r, b^\ast_s\}$. This flat must be a connected basis flat for $b$ by Lemma \ref{typeofflats}, hence there exists some element $b_j\in b$ satisfying $\cl\{b_i, b_j\} = \cl\{b_k^\ast, b_l^\ast\} $. 
If neither $b_i$ nor $b_j$ are contained in $b^\ast$, then these two vertices are linked with two edges in $G_b(b^\ast)$. Therefore $b \cap b^\ast = \emptyset$ implies that every vertex in the connected graph $G_b(b^\ast) $ has some neighboring vertex to which it is linked by at least two edges. Since $|b| = |b^\ast| = |\rk(M) | \geq 3$, this is impossible. 

Therefore we see that $b \cap b^\ast$ is a singleton set $\{b_0\}$, and there are $d = \rk(M) -1$ edges in $G_b(b^\ast)$.
Since $G_b(b^\ast)$ is a connected graph with $d+1$ vertices, it can not contain any double edge. We have seen previously that every $b_i \in b \backslash \{b_0\}$ is contained in some rank $2$ flat of the form $\cl\{b_k^\ast, b_l^\ast\}$, which has to equal to $\cl\{b_i, b_0\}$. Therefore $b_i$ is connected by an edge to $b_0$, which implies that $G_b(b^\ast)$ has the shape of a star with center $b_0$.

It remains to consider the case that $H$ is a proper connected component of $G_b(b^\ast)$.
By our initial observation we know that $H = G_b(I)$ for a proper subset $I$ of $b^\ast$ that contains at most one element of $b$. Since $M$ is connected, there exists an element $b_i \in V(H)= \supp_b(I)$ and an element $b_j \in b \backslash V(H)$ such that $F_{ij} \neq \emptyset$. Let $e \in F_{ij}$. Note that $e \notin b^\ast$. Since $b^\ast$ is a Cremona basis, we have $e \in \cl\{b_k^\ast, b_l^\ast\}$ for some $k,l$.
Assume that this flat is a non-basis flat for the basis $b$. Then by Lemma \ref{nonbasisflats} (3) its support graph with respect to $b$ is a triangle containing the edge given by $e$ between the vertices $b_i$ and $b_j$. In this case neither $b_i$ nor $b_j$ are contained in $\cl\{b_k^\ast, b_l^\ast\}$, so that $b_k^\ast$ and $b_l^\ast $ must give the other two edges in the triangle. This is impossible, since $b_i$ and $b_j$ are supposed to belong to different connected components of $G_b(b^\ast)$. Hence  we find that $\cl\{b_k^\ast, b_l^\ast\}$ is a basis flat for the basis $b$. Since its support contains $b_i$ and $b_j$, we find that the rank two flats $\cl\{b_i,b_j\}$ and $\cl\{b_k^\ast, b_l^\ast\}$ coincide. Since $b_i$ and $b_j$ are supposed to belong to different connected components of $G_b(b^\ast)$, neither $b_k^\ast$ nor $b_l^\ast$ can lie in $b^\ast \backslash b$. This implies that $b_i$ lies in $b \cap b^\ast$, so that $V(H) \cap b^\ast = \{b_i\}$. Let $b_k \neq b_i$ be an element of $V(H) = \supp_b(I)$.  Then $b_k \in b \backslash b^\ast$, and it is connected to some $b_l \in V(H)$ by an edge in $b^\ast$. The same argument as above shows that there is a double edge between $b_k$ and $b_l$, if $b_l \neq b_i$. Since we have already seen that $|b \cap b^\ast|$ is equal to the number of connected components of $G_b(b^\ast)$, which is a graph with $\rk(M)$ many vertices and $\rk(M) - |b \cap b^\ast|$ many edges, this is a contradiction. Hence $H$ is simple and star-shaped with center $b_i$.

 Moreover, this argument also shows that if $b_i$ and $b_j$ are vertices in different components of $G_b(b^\ast)$ such that $F_{ij}\neq \emptyset$, we find that $b_i$ and $b_j$ are in fact contained in $b \cap b^\ast$.
 \end{proof}

Any matroid automorphism maps a Cremona basis to a Cremona basis. The next result shows that this action of $\Aut(M)$ on the set of Cremona bases is transitive. 
\begin{corollary}\label{cor: involution}
   Assume that the simple and connected matroid $M$ of rank $d+1 \geq 3$ admits two different Cremona bases $b$ and $b^\ast$. Then there exists an involutive matroid automorphism $\varphi\colon E(M) \rightarrow E(M)$ mapping $b$ to $b^\ast$.
\end{corollary}
\begin{proof} 
By Theorem \ref{twocremonas} every element $c \in b \backslash b^\ast$ is connected to a unique vertex in $b \cap b^\ast $ in the support graph $G_b(b^\ast)$. We denote by $e_b(c) \in b^\ast \backslash b$  the label of the corresponding edge in $G_b(b^\ast)$. The element $e_{b^\ast}(c)$ for $c \in b^\ast \backslash b$ is defined in an analogous way. Then we define the map $\varphi\colon E(M) \rightarrow E(M)$ by 
\begin{align*}
      \varphi(c) & =  c & \mbox{for all } c \in b \cap b^\ast \\
   \varphi(c)  & =  c & \mbox{for all }
 c \in E(M) \backslash (b \cup b^\ast) \\
 \varphi(c) & = e_b(c) & \mbox{ if }c \in b \backslash b^\ast \\
 \varphi(c)&  = e_{b^\ast}(c) & \mbox{ if } c \in b^\ast \backslash b.  
 \end{align*}
Note that for every $c \in b \backslash b^\ast$ with associated edge element $e_b(c) \in b^\ast \backslash b$ there is an element $b_0 \in b \cap b^\ast$ such that $e_b(c) \in \cl\{b_0, c\}$,  which implies that  $\{ b_0, c, e_b(c) \}$ is a circuit. Thus $c$ is the edge in $G_{b^\ast}(b)$ between $b_0$ and $e_b(c)$, hence $c = e_{b^\ast}(e_b(c))$. Analogously, $c = e_b(e_{b^\ast}(c))$ for every $c \in b^\ast \backslash b$.


Hence $\varphi$ is indeed an involution as a map of sets and hence bijective. 
Since $\varphi(b) \subset b^\ast$ by definition, we deduce that $\varphi$ restricts to a bijection between $b$ and $b^\ast $.

Let us now show that for any two different elements $b_i, b_j$ of the basis $b$ we have $\varphi (\cl\{b_i,b_j\} ) = \cl\{\varphi(b_i), \varphi(b_j)\}$. This is obvious by definition of $\varphi$ if both $b_i$ and $b_j$ are contained in $b \cap b^\ast$. Hence let us assume that least one of them, say $b_i$,  is contained in $b \backslash b^\ast$.
By Theorem \ref{twocremonas}, we find that $F_{ij} = \emptyset$, if $b_i$ and $b_j$ are vertices in different components of the support graph $G_b(b^\ast)$. By definition, $\varphi(b_i)$ is an edge in the same component of $G_b(b^\ast)$ as $b_i$. A parallel argument applied to the elements $\varphi(b_i)$ and $\varphi(b_j)$ of $b^\ast$ shows that in this case $\cl\{\varphi(b_i), \varphi(b_j)\} = \{\varphi(b_i), \varphi(b_j)\} $,
which is equal to 
$\varphi(\cl\{b_i,b_j\})$.

Now assume that $b_i \in b \backslash b^\ast$ and $b_j$ are vertices in the same component of $G_b(b^\ast)$, so that we either have $b_j\in b \cap b^\ast$ or $b_j \in b \backslash b^\ast$. In the first case we find $\varphi(b_i) \in \cl\{b_i,b_j\}$, which implies $\cl\{b_i, b_j\} = \cl\{\varphi(b_i), b_j\} = \cl\{\varphi(b_i), \varphi(b_j)\}$. 

Hence we may assume that $b_i$ and $b_j$ are both contained in $b \backslash b^\ast$. Let us first treat the case that $F_{ij} \neq \emptyset$, i.e.~that the flat $\cl\{b_i, b_j\}$ is connected. By Theorem \ref{twocremonas}, the connected rank $2$ flat $\cl\{b_i, b_j\}$ cannot contain two elements of $b^\ast$ as this would lead to a double edge in the support graph $G_b(b^\ast)$. Therefore it is a non-basis flat for the Cremona basis $b^\ast$, and its support graph for $b^\ast$ is a triangle by Lemma \ref{nonbasisflats} (3). Hence there can be only one element $e \in F_{ij}$ if this set is non-empty. 
We claim that $e \in \cl\{\varphi(b_i), \varphi(b_j)\}$. Indeed, since $b^\ast$ is also a Cremona basis, there exist $b_k, b_l \in b \backslash b^\ast$ such that $\{ \varphi(b_k), \varphi(b_l), e \}$ is a circuit, which implies that $\{ b_0, b_k, b_l, e \}$ is dependent. Using that $b$ is a Cremona basis, this is possible only if $\{ b_k, b_l \} = \{ b_i, b_j \}$. Thus the flat $\cl\{\varphi(b_i), \varphi(b_j)\}$ is connected and its support graph with respect to $b$ has 3 vertices. By Lemma \ref{nonbasisflats} (3), its support graph is a triangle and hence $\cl\{\varphi(b_i),\varphi(b_j)\} = \{\varphi(b_i), \varphi(b_j), e\} = \varphi(\cl\{b_i,b_j\})$.

A parallel argument shows that connectedness of $\cl\{\varphi(b_i), \varphi(b_j)\}$ implies connectedness of $\cl\{b_i,b_j\}$. Hence our claim also follows for $F_{ij} = \emptyset$.




Now consider a connected flat $F$ which is a basis flat for the Cremona basis $b$. Then $F$ is generated by $F \cap b$, which is mapped bijectively to the set $\varphi(F) \cap b^\ast$ via $\varphi$.

Hence we find that 
\begin{eqnarray*} \varphi(F) & = & \varphi(\cl(F \cap b)) = \varphi\bigg(\bigcup_{b_i, b_j \in F \cap b} \cl\{b_i,b_j\}\bigg)
=\bigcup_{b_i, b_j \in F \cap b} \varphi(\cl\{b_i,b_j\})\\ & =& \bigcup_{b_i,b_j \in F \cap b} \cl\{\varphi(b_i), \varphi(b_j)\} = \bigcup_{b_i^\ast, b_j^\ast \in \varphi(F) \cap b^\ast} \cl\{b_i^\ast, b_j^\ast\} = \cl(\varphi(F) \cap b^\ast),
\end{eqnarray*}
which proves that the set $\varphi(F)$ is in fact a flat
of rank $ |\varphi(F) \cap b^\ast| = \rk(F)$. 


Connected non-basis flats disjoint from $b^\ast$ will be mapped to themselves under $\varphi$, hence their rank is preserved under $\varphi$.

Now consider a connected non-basis flat $F$ for the Cremona basis $b$ containing an element of $b^\ast$. Hence $F$ is a basis flat for the Cremona basis $b^\ast$, and the same argument as above exchanging the roles of $b$ and $b^\ast$ shows that $\varphi(F) $ is a flat of rank  $\rk(\varphi(F)) = |\varphi(F) \cap b| = |F \cap b^\ast| = \rk(F)$.

Any flat $F$ in $M$ is a disjoint union of connected flats $F_i$. Since $\varphi$ is bijective, the set $\varphi(F)$ is the disjoint union of the flats $\varphi(F_i)$. Therefore we deduce $\rk(\varphi(F)) \leq \sum_i \rk(\varphi(F_i)) = \sum_i \rk(F_i) = \rk (F)$. 

We can apply this to the flat $F' = \cl(\varphi(F))$ and find $\rk(\varphi(F')) \leq \rk(F')$.
Since $\varphi$ is an involution, we have $\rk(F) \leq \rk(\varphi(\cl(\varphi(F)))) = \rk(\varphi(F'))$,
which implies $ \rk(F) \leq \rk(F') = \rk(\varphi(F))$.
Hence we showed that for every flat $F$ its image $\varphi(F)$ is a flat of the same rank. Therefore  $\varphi$ is a rank-preserving bijection on the set of flats of $M$. 
\end{proof}

We will now prove that matroids admitting two different Cremona bases are realizable over large enough fields. We need the following notation: If $M$ admits the Cremona basis $b = \{b_0, \ldots, b_d\}$ with first basis vector $b_0$, we define $E_+ = \cup_{0 <i <j} F_{ij}$ and $E_0 = \cup_{i \neq 0} F_{0i}$. Then $E(M) = b \cup E_0 \cup E_+$. 
We define an equivalence relation on $E_0$ by \[e \sim f, \mbox{ if and only if } E_+ \vee\{e\} = E_+ \vee \{f\}.\] 

\begin{thm}\label{thm:repcrit} Let $M$ be simple and connected of rank $d+1 \geq 3$ and assume that $M$ admits two different Cremona bases $b = \{b_0, \ldots, b_d\}$ and $b^\ast = \{ b_0, b_1^\ast, \ldots, b_d^\ast\}$ such that $b \cap b^\ast$ consists of the element $b_0$. 

Let $N $ be the number of equivalence classes in $E_0$ for the equivalence relation defined above. Then  $M$ is realizable over all fields with at least $N+1 $ elements, in particular over all fields with at least $|E(M)| - d$ elements.

More precisely, we prove that $M$ can be realized over any sufficiently large field $K$ by a subset of the vectors $\{ v_i \mid 0 \le i \le d \} \cup \{ v_0 - av_i \mid a \in K, 1 \le i \le d \} \cup \{ v_i - v_j \mid 1 \le i < j \le d \}$, where $v_0, \dotsc, v_d$ is a basis of $K^{d+1}$.
\end{thm}
\begin{proof} Let $K$ be any field with at least $N+1$ elements. Note that  $N \leq |E_0| \leq |E(M)| - d - 1$, so that a field with at least $|E(M)| - d$ elements satisfies this condition.

(1) To start the proof, we have a look at $E_+ = \bigcup_{0<i<j} F_{ij}  = \cl(b\backslash \{b_0\}) \backslash b$. By Theorem \ref{twocremonas}, every element in $b^\ast \backslash \{ b_0 \}$ is contained in $\cl\{ b_0, b_i \}$ for some $b_i \in b$ and thus $E_+ \cap b^\ast = \emptyset$. Since $b^\ast$ is also a Cremona basis, we find that every element $e$ in $E_+$ is contained in the closure of two elements of the basis $b^\ast$. Note that neither of these elements may be $b_0$, since otherwise the closure of $e$ and $b_0$ would be a connected basis flat for the basis $b$, which contradicts the choice of $e$. This implies that $E_+ \subset \cl(b \backslash \{b_0\}) \cap \cl(b^\ast \backslash \{b_0\})$. This inclusion is in fact an equality, since
$\cl(b^\ast \backslash \{ b_0 \}) \cap b = \emptyset$ as above by Theorem \ref{twocremonas}.
Hence $E_+$ is a flat. 

Now let $e \in F_{0i}$ be any element in $E_0$, and assume that $b_i$ is in the support of $E_+$. Let $a \subset b$ be the vertices in the connected component of the support graph $G_b(E_+)$ containing the vertex $b_i$. Then $D =  \cl(a) \backslash a = \cl(a) \cap E_+$ has the property that the support graph of $D \vee \{e\}$ is connected.  If $D \vee \{e\}$ was a basis flat for the Cremona basis $b$, its rank would be equal to $|\supp_b(D \vee \{e\})| = |\supp_b(D)| +1$. But $D$ is a non-basis flat for $b$, and hence Lemma \ref{nonbasisflats} (4) gives a contradiction. Therefore $D \vee \{e\}$ must be a non-basis flat for $b$. 

(2) We choose an injective map $\kappa$ from the set of equivalence classes $E_0 /{\sim}$ to $ K^\times$, and denote the image of the equivalence class of some $e \in E_0$ simply by $\kappa(e)$.  We denote by $v_0, \ldots, v_d$ be canonical basis of $K^{d+1}$.  Now we define an embedding $\sigma$ of $E(M)$ into $K^{d+1}$  as follows: 
\begin{align*}
\sigma(b_i) & = v_i & \mbox{for all } b_i \in b \\
\sigma (e) & = v_0 - \kappa(e) v_i & \mbox{if }e \in F_{0i} \mbox{ for some } i \in \{1, \ldots d\}\\
\sigma(e) & = v_i - v_j & \mbox{for all }e \in F_{ij} \mbox{ with }0 < i < j.
\end{align*}

(3) Let us first check that $\sigma$ is indeed injective. This is clear on $b$. Let $0<i<j$ be indices such that $F_{ij} \neq \emptyset$. By Theorem \ref{twocremonas},  
the connected rank $2$ flat $\cl\{b_i, b_j\}$ cannot contain two elements of $b^\ast$ as this would lead to a double edge in the support graph $G_b(b^\ast)$. Therefore it is a non-basis flat for the Cremona basis $b^\ast$, and its support graph for $b^\ast$ is a triangle by Lemma \ref{nonbasisflats} (3). Hence there can be only one element in $F_{ij}$. 
This shows that the map $\sigma$ is also injective on $E_+$. It remains to check injectivity on $E_0$. Consider equivalent elements $e,f$ in $F_{0i}$ for some $i \geq 1$. Then $E_+ \vee \{e\}$ equals $E_+ \vee \{f\}$. Assume that $e \neq f$, which implies $\cl\{e,f\} = \cl\{b_0,b_i\}$. Then there is a  support-connected subflat of $E_+ \vee \{e\}$  containing $e$ and $f$. This must be a basis flat for $b$, which is impossible by (1).

(4) It remains to show that $\sigma$ is in fact an isomorphism between $M$ and the realizable matroid $M'$ given by the image $\sigma(E(M))$. Note that $c := \{ v_0, \ldots, v_d \}$ is a Cremona basis of the matroid $M'$, and that by definition
\[\sigma(\cl\{b_i,b_j\}) = \cl\{v_i, v_j\} \mbox{ for all } i \neq j.\]
This implies that $\sigma$ maps the support graph of a subset $F$ of $M$ with respect to $b$ isomorphically to the support graph of the subset $\sigma(F)$ of $M'$ with respect to $c$.

(5) Assume that $F$ is a connected  non-basis flat for $M$ containing $b_0$ in its support. 
Let $f_i \in F \cap F_{0i}$ and $f_j \in F \cap F_{0j}$ for some $0 < i < j$. 
Then $b_i$ and $b_j$ lie in the same connected component of $F \cap E_+$. Indeed,
let $d \subset b$ be the set of vertices in the connected component of $G_b((F \cap E_+) \cup b)$ containing $b_i$, then $F = (F \cap \cl(d \cup \{ b_0 \})) \cup (F \cap \cl(b \setminus d))$ is the disjoint union of two flats. As $F$ is connected, we find that  $b_j$ has to be  contained in the same connected component of $G_b((F \cap E_+) \cup b)$ as $b_i$. Thus $F \cap E_+$ is support-connected, so that $\rk(F) = \rk(F \cap E_+) + 1$ by Lemma \ref{nonbasisflats}.  We deduce $(F \cap E_+) \vee \{ f_i \} = (F \cap E_+) \vee \{ f_j \}$ and hence $f_i \sim f_j$.
Hence $a := \kappa(f_i)$ is independent of the choice of  $f_i \in F \cap E_0$. 

(6) We show that  $\rk(\sigma(F)) \le \rk(F)$ for every connected flat $F$ of $M$. Since $\sigma$ respects the supports, this is clear if $F$ is a basis flat for $b$. If $F$ is a non-basis flat for $b$ with $b_0 \notin \supp_b(F)$, then $\sigma(F)$ is contained in  $\langle v_i - v_j \mid b_i, b_j \in \supp_b(F) \rangle$, which is a subspace of dimension $|\supp_b(F)| - 1 = \rk(F)$. On the other hand, if $b_0 \in \supp_b(F)$, then by (5) $a = \kappa(f)$ is independent of the choice of $f \in F \cap E_0$ and $\sigma(F)$ is contained in $\langle v_0 - av_i, v_i - v_j \mid b_i, b_j \in \supp_b(F) \setminus \{ b_0 \} \rangle$, which is a subspace of dimension $|\supp_b(F)| - 1 = \rk(F)$.


(7) We show that $\rk(\sigma^{-1}(F')) \le \rk(F')$ for every connected flat $F'$ of $M'$. Since $\sigma^{-1}$ respects the supports, this is clear if $F'$ is a basis flat for $c$. If $F'$ is a non-basis flat for $c$ with $v_0 \notin \supp_c(F')$, then $\sigma^{-1}(F')$ is contained in the non-basis flat $E_+ \cap \sigma^{-1} \cl(\supp_c(F'))$, which has rank $|\supp_c(F')|-1 = \rk(F')$. On the other hand, if $v_0 \in \supp_c(F')$, the subflat $F'_1 := F' \cap \cl(c \setminus \{ v_0 \})$ is a support-connected non-basis flat such that $F' = F'_1 \vee \{ f' \}$ for all $f' \in F' \setminus F'_1$. This implies that for every two elements $v_0 - av_i$ and $v_0 - bv_j$ in $F' \setminus F'_1$ we have $a = b$, hence $e \sim e'$ for all $e, e' \in \sigma^{-1}(F') \cap E_0$. Thus for any $e \in \sigma^{-1}(F') \cap E_0$ we find that $\sigma^{-1}(F')$ is contained in the flat $(E_+ \vee \{ e \}) \cap \sigma^{-1} \cl(\supp_c(F'))$, which is a non-basis flat by (1) and thus has rank $|\supp_c(F')| - 1 = \rk(F')$.

(8) Any flat $F$ in $M$ is a disjoint union of connected flats $F_i$. Using (6) we find $\rk(\sigma(F)) \leq \sum_i \rk(\sigma(F_i)) \leq \sum_i \rk (F_i) = \rk(F)$. In a similar fashion we show that $\rk(\sigma^{-1}(F')) \leq \rk(F')$ for all flats $F' $ in $M'$. Hence we find $\rk(\sigma^{-1}(\cl(\sigma(F)))) \leq \rk(\sigma(F)) \leq \rk (F)$. As $F$ is contained in $\sigma^{-1}(\cl(\sigma(F)))$, we have equality everywhere, which implies that  $\sigma(F)$ is a flat of the same rank as $F$.   From this we deduce that $\sigma$ is an isomorphism of matroids. 
\end{proof}

The following example shows that the hypothesis $|b \cap b^\ast| = 1$ of the preceding theorem is necessary for the desired realizability statement. 

\begin{example}
The Dowling matroid $Q_3(\Z/2\Z \times \Z/2\Z)$ of rank $3$ given by the group $G = \Z/2\Z \times \Z/2\Z$ admits a unique Cremona basis $\{ b_1, b_2, b_3 \}$ and is not realizable over any field by \cite[Theorem 6.10.10]{ox}. Consider the parallel connection $M$ of $Q_3(\Z/2\Z \times \Z/2\Z)$ and the uniform matroid $U_{2,3}$ on the ground set $\{ 1, 2, 3 \}$ obtained by identifying $b_1$ and $1$. Then $\{ b_1, b_2, b_3, 2 \}$ and $\{ b_1, b_2, b_3, 3 \}$ are two different Cremona bases of $M$ whose intersection has more than one element.
\end{example}

\section{Coxeter arrangement matroids}
\label{sec:Cox}
In this section, we study matroids given by Coxeter arrangements. To be precise, let $(W,S)$ be a Coxeter system with finite Coxeter group $W$, which is generated by the set $S$ consisting of elements of order $2$ satisfying $(s s')^{m(s,s')} = 1$,  see \cite[Section 1.9]{hu}. Let us recall some facts about the geometric representation of $W$.  Denote by $(\alpha_s)_s$ the canonical basis of $V= \mathbb{R}^S$. Then $B(\alpha_s, \alpha_{s'}) = -\cos(\pi / m(s,s'))$ defines a positive definite symmetric bilinear form on $V$, see \cite[Section 6.4, Theorem]{hu}. There exists a unique homomorphism $\sigma\colon W \rightarrow GL(V)$ mapping $s$ to the reflection $\sigma(s) (x) = x - 2 B(\alpha_s, x) \alpha_s$. It identifies $W$ with a subgroup of the orthogonal group with respect to $B$, see \cite[Section 5.3]{hu}. 

We assume that $W$ is irreducible, so that the Coxeter graph associated to $(W,S)$ is isomorphic to one of the types listed in \cite[Section 2.4]{hu}. If $W$ is crystallographic, i.e. if $\sigma(W)$ stabilizes a lattice in $V$, then  there exists a root system $\Phi$ with associated Weyl group $W$, see \cite[Section 2.9]{hu}. For an introduction to root systems see \cite[Chapter VI]{bou}. 

Consider the finite set $\Phi(W,S) = \{ \sigma(w)(\alpha_s): w \in W, s \in S\} \subset V$. We also write $w(\alpha_s) = \sigma(w) (\alpha_s)$. We put 
\[\Phi^+(W,S) = \{v \in \Phi(W,S): v = \sum_{s \in S} c_s \alpha_s: c_s \geq 0 \mbox{ for all } s \in S\}.\]
The matroid $M = M(W,S)$ associated to the finite Coxeter system $(W,S)$ is defined as the realizable matroid given by the subset $\Phi^+(W,S)$ of $V$. Note that scaling the elements in $\Phi^+(W,S)$ gives another realization of $M$, which is sometimes convenient for explicit calculations. The Weyl group $W$ acts on $M$ via
\[ w(v) = \left\{\begin{array}{rl}
  w(v)   & \mbox{if }  w(v) \in \Phi^+(W,S)\\
  -w(v)   & \mbox{if } {-w(v)} \in \Phi^+(W,S) 
\end{array}
\right.\]



We have seen in section \ref{sec:A} that the coarse Bergman fan associated to the Coxeter arrangement of type $A_n$ admits several combinatorial Cremonas, and that its automorphism group is generated by matroid automorphisms and one Cremona map.

Let us now discuss the case of a crystallographic Coxeter group of type $B_n$. The associated realizable matroid is given by 
\[E(M) = \Phi^+ = \{x_i \mid 1 \le i \le n \} \cup \{ x_i \pm x_j \mid 1 \le i < j \le n \}
 \]
in $\mathbb{R}^n$, where $x_i$ denotes the $i$-th canonical basis vector. Then $b = \{x_1, \ldots, x_n\}$ is a Cremona basis for $M$ in the sense of \cite[Definition 8.1]{sw}, i.e.~$\crem_b$ is a linear automorphism of $B_c(M)$. 


\begin{prop}\label{prop:B}
Let $M$ be the matroid associated to the root system of type $B_n$ for $n\geq 3$. Let $\varphi\colon B_c(M) \rightarrow B_c(M)$ be a linear automorphism of the Bergman fan. Then the automorphism group
$\Aut(B_c(M))$ is isomorphic to $(\Z/2\Z)^{n} \rtimes S_n$.
\end{prop}
\begin{proof}

Note that $M$ is isomorphic to the Dowling matroid on $n$ elements for the group $G = \mathbb{Z} / 2 \mathbb{Z}$, see \cite[Lemma 6.10.11]{ox} applied to $\{ -1, 1 \}$ as a finite subgroup of $\mathbb{R}^\times$.  Hence, using  \cite[Proposition 10.3]{sw} if $n \geq 4$, and \cite[Theorem 9.5]{sw} if $n =3$, we find that $\Aut(B_c(M))$ is generated by $\Aut(M)$ and the Cremona map $\crem_b$.

Note that the support graph $G_b(M)$ has two edges $x_i \pm x_j$ between any two vertices $x_i$ and $x_j$. Assume that $b^\ast$ is another Cremona basis.  It is easy to see that $G_b(b^\ast)$ is connected. Applying Theorem \ref{twocremonas} we find that there are at least two elements $x_k, x_l \in b\backslash b^\ast$. They generate a connected non-coordinate flat for the basis $b^\ast$, which by Lemma \ref{nonbasisflats} cannot contain four elements.
Therefore we find that $b = \{x_1, \ldots, x_n\}$ is the only Cremona basis for $M$.

Since a matroid automorphism $\varphi$ of $M$ maps a Cremona basis to a Cremona basis, it satisfies $\varphi(b) = b$. This implies that $\crem_b \circ \varphi = \varphi \circ \crem_b$ in the automorphism group of $B_c(M)$. 

 As the Cremona map has order two, we find that $\Aut(B_c(M))$ is isomorphic to the direct product $\Aut(M) \times \Z/2\Z$.

 By \cite[Theorem 1.2]{sft} we have $\Aut(M) \cong (\Z/2\Z)^{n-1} \rtimes S_n$, which implies our claim.

\end{proof}

We will now show that for Coxeter arrangements beyond $A_n$ and $B_n$ every automorphism of the coarse Bergman fan of the associated matroid is given by a matroid automorphism. The first observation is that the coarse structure on the Bergman fan of $M = M(W,S) $
coincides with the minimal nested set structure $B_m(M)$ by \cite[Theorem 1.2]{arw}. The rays in the fan structure of $B_m(M)$, and hence also the rays in coarser fan structures, are given by connected flats $F$ of $M$. The rank or corank of a ray denotes the rank or corank of the associated connected flat, where the corank of $F$ is defined as $\rk(M) - \rk(F)$. 

We need the following Lemma from \cite{arw}:

\begin{lemma}\label{flats} Let $M$ be the matroid associated to a finite irreducible Coxeter system $(W,S)$ with Coxeter graph $\Gamma$, and let $F$ be a flat of $M$. Then the restriction $M|_{F}$ is isomorphic to the  matroid associated to the Coxeter system given by a subset $J$ of $\Gamma$. The flat $F$ is connected if and only if $J$ spans a connected subgraph of $\Gamma$.

\end{lemma}
\begin{proof} By, \cite[Section 5]{arw},  every flat $F$ of $M$ is of the form $F = w \Phi^+_J$ for some $w \in W$ and a subset $J \subset S$. Here $\Phi^+_J \subset \Phi^+ \subset V$ denotes the intersection of $\Phi^+$ with the span of $J$.
Obviously $F$ is connected if and only if $J$ spans a connected subgraph of $\Gamma$.
\end{proof}

\begin{lemma}\label{corank} Let $M$ be the matroid associated to an irreducible finite Coxeter system $(W,S)$, and let $\varphi\colon B_c(M) \rightarrow B_c(M)$ be an automorphism of the Bergman fan. Then $\varphi$ maps rank one rays of $B_c(M)$ either to rank one or to corank one rays. 

\end{lemma}

\begin{proof} We may assume that the rank of $M$ is at least $4$, otherwise the statement is trivial. Our claim follows from \cite[Corollary 7.5]{sw}, once we show that for every element $i$ of the ground set $E = E(M)$ the matroid $M/i$ is not a non-trivial parallel connection. Hence let us assume that $M/i$ is a parallel connection of $N_1$ and $N_2$ along some joint element $j$, where $N_1$ and $N_2$ have ranks at least $2$. Then $F = \cl\{i,j\}$ is a flat of rank $2$ in $M$.

By Lemma \ref{flats} we find that $F = w^{-1} \Phi_J^+$ for some $w\in W$ and $J \subset S$. Hence after applying $w$ (which induces an automorphism of $M$ and hence of $B_c(M)$), we may assume that $F = \Phi_J^+$ for a subset $J$ of $S$  of cardinality  $2$. By definition of $\Phi_J^+$ there exists an element $w \in W_J$ and an element $s \in J$ such that $i = w(\alpha_s)$. After applying $w^{-1}$ (which maps $F$ to itself) we  may therefore assume that $i\in S$. Here we abuse notation to identify any $s \in S$ with the corresponding element $\alpha_s \in E(M)$.

By  \cite[Lemma 5.3]{arw} there exists an element $e \in E(M)$ such that $S \cup \{e\}$ is a circuit in $M$ containing $i$.
Then $S\backslash\{i\} \cup \{e\}$ is a circuit in $M/i$. This circuit has full rank, so it cannot be contained in $N_1$ or $N_2$. Hence it is made up from circuits in both of them by deleting $j$. In particular, $e \neq j$. Without loss of generality we have  $e \in N_1 \backslash N_2$. Then we find $S_1$ and $S_2$ such that $S_1 \cup S_2 \cup \{i\} =S$ is a disjoint union, and such that $S_1 \cup \{e, j\}$ is a circuit in $N_1$ and $S_2 \cup \{j\}$ is a circuit in  $N_2$ and hence in $M/i$.
Hence $S_2 \cup \{i,j\}$ is a circuit in $M$ or $S_2 \cup \{j\}$ is a circuit in $M$.
Looking at the linear combination of $j$ with respect to the basis $(\alpha_s)_{s \in S}$, we find that either all vectors given by $S_2 \cup \{i\}$ or all vectors given by $S_2$ contribute. On the other hand, $j$ is contained in $F$, so that we can write it as a linear combination of the two elements given by $J$,  one of which is $i$.
This implies that the cardinality of $S_2 $ is one. As $S_1 \cup \{e, j\}$ is a circuit in $N_1$, we find that $\rk(N_1) \geq |S_1|+1 = |S| -1 = \rk(M) -1 = \rk(M/i)$.
This contradicts the fact that $M/i$ is a non-trivial parallel connection of $N_1$ and $N_2$. 
\end{proof}

\begin{thm}\label{auto}
Let $M$ be the matroid associated to a Coxeter arrangement which is not of type $A_n$ or $B_n$.
  Then every automorphism $\varphi\colon B_c(M) \rightarrow B_c(M)$ is induced by a matroid automorphism of $M$. 
\end{thm}

Note that the matroid automorphisms of $M$ are determined in \cite{sft}.
\begin{proof} 
For rank $2$ matroids this is always true. Hence we may assume that  $M$ has rank at least $3$. 

Recall that by \cite[Theorem 1.2]{arw}, the coarse structure on the Bergman fan coincides with the minimal nested set structure.

We look at the subgraph $S$ of the $1$-skeleton of $B_c(M)$ given by all rays associated to rank one or corank one flats. Hence the vertices of $S$ are given by the connected rank $1$ and corank one flats in $M$, and two vertices are connected by an edge, if and only if the associated rays form a cone in the minimal nested set structure.  Hence the neighbors of a rank one ray $F$ in $S$ are given by all connected corank one rays containing $F$ plus all rank one rays $G$ such that $F \vee G$ is disconnected. The neighbors of a corank one ray $F$ on the other hand are given by the elements of $F$, since $M$ is simple.  By Lemma \ref{corank},  every automorphism $\varphi\colon B_c(M) \rightarrow B_c(M)$ induces an automorphism $\varphi_S$ of the graph $S$. 

Note that it suffices to show that every rank one vertex in $S$ has strictly more neighbors  than any corank one vertex in $S$. Indeed, if this is the case, then  $\varphi_S $ has to map every vertex given by a rank one ray to a vertex given by a rank one ray, which implies that $\varphi$ preserves the set of rank one rays in the minimal nested set structure of the Bergman fan. This provides a map $E(M) \rightarrow E(M)$, which has to be a matroid automorphism as the linear map $\varphi$ maps indicator vectors of flats to indicator vectors of flats (modulo $\R \mathbbm{1}$). This matroid automorphism induces $\varphi$ on the Bergman fans.

By the classification of Coxeter systems of rank at least three (see e.g. \cite[Section 2.4]{hu}), the remaining possible types are  $D_n$ for $n \geq 3$, $E_6$, $E_7$, $E_8$, $F_4$, $H_3$ and $H_4$, as $C_n$ gives the same matroid as $B_n$. We will use the realizations for root systems given in 
\cite[VI, \S 4]{bou} and the realizations for type $H$ in \cite{H3} and \cite{H4} to discuss these cases. 

We will show in a case-by-case analysis that every rank $1$ vertex in $S$ has strictly more neighbors  than any corank $1$ vertex in $S$.

If $M$ is associated to a root system $\Phi$ of type $D_n$ for $n \geq 3$, then $E(M)$ can be realized as the set $\{x_i \pm x_j :  i <j\}$ with respect to the canonical basis $x_1, \ldots, x_n$ of $\mathbb{R}^n$. We begin by listing the neighbors of the flat given by the element $x_1 + x_2$. Among those are the rank one flats given by the elements $x_i \pm x_j$ with $2 <i < j$ plus the element $x_1 - x_2$, which are altogether $(n-2)(n-3)+1$ many flats. We also find $n-2$ corank one neighbors of $x_1 + x_2$ of the form 
\[\{ x_i \pm x_j \mid i<j \mbox{ and } i, j \ne 3 \}, \dotsc, \{ x_i \pm x_j \mid i<j \mbox{ and } i, j \ne n \}\](which are all isomorphic to matroids of type $D_{n-1}$), and $2^{n-2}$ corank one neighbors of the form $\cl\{ x_1 + x_2, x_1 + \epsilon_3x_3, \dotsc, x_1 + \epsilon_nx_n \}$ for choices of signs $\epsilon_3, \ldots, \epsilon_n$ (which are all isomorphic to matroids of type $A_{n-1}$). As the Weyl group of $\Phi$ acts simply transitively on $M$ by \cite[VI, \S 1 no. 3, Prop. 11]{bou}, we find that every vertex of $S$ given by a rank 1 ray has at least $(n-2)(n-3)+1  + (n-2) + 2^{n-2} >(n-2)^2 + 2^{n-2}$ neighbors in $S$. On the other hand, any connected flat $H$ in $M$ of corank one gives rise to a root system such that the Coxeter graph is a connected subgraph of the $D_n$-graph by Lemma \ref{flats}, hence it contains at most $(n-1)(n-2)$ elements. A quick estimate shows that this is smaller or equal to  $(n-2)^2 + 2^{n-2}$, which proves our claim for $D_n$.

The matroid associated to the root system of type $F_4$ can be realized in $\mathbb{R}^4$ by the vectors $x_1, \dotsc, x_4$ plus 12 vectors of the form $x_i \pm x_j$ for $i < j$ plus 8 vectors of the form $x_1 + \epsilon_2x_2 + \epsilon_3x_3 + \epsilon_4x_4$ with signs $\epsilon_2, \epsilon_3, \epsilon_4 \in \{ -1, 1 \}$. Every connected hyperplane of $F_4$ is isomorphic to the matroid associated to $B_3$ and thus contains 9 elements. On the other hand, every rank one ray has degree 15 in the graph $S$.

Let us now consider the matroid associated to the root system of type $E_8$, which can be realized in $\mathbb{R}^8$ as the $56$ vectors of the form $x_i \pm x_j$ ($1 \le i < j \le 8$) plus the $64$ vectors of the form $x_1 + \epsilon_2x_2 + \dotsb + \epsilon_8x_8$ with signs $\epsilon_2, \dotsc, \epsilon_8 \in \{ -1, 1 \}$ such that $\prod_{i=2}^8 \epsilon_i =1$ by \cite[VI, \S 4 no. 10]{bou}.   Using Lemma \ref{flats}, we find that the maximal cardinality of a connected hyperplane in $E_8$ is the cardinality of the matroid induced by $E_7$, i.e.  $63$. On the other hand, every rank one ray in $S$ has exactly 63 rank one neighbors: As the Weyl group acts transitively on $E(M)$,  it suffices to count the rank 1 neighbors of the root $x_1 + x_2$, which are $x_1 - x_2$, 30 roots of the form $x_i \pm x_j$ with $3 \leq i < j \leq 8 $  and 32 roots of the form $x_1 + \epsilon_2x_2 + \dotsb + \epsilon_8x_8$ with $\epsilon_2 = -1$. But every element also has corank one neighbors, hence we conclude by the usual argument. 

For $M$ of type $E_7$, a similar argument shows that every rank one vertex in $S$ has $30$ rank one neighbors, 
and the maximal cardinality of a connected hyperplane is $36$, so it remains to show that every rank one vertex in $S$ has at least 7 corank one neighbors. Consider the subgroup $H := \tau S_6 \tau^{-1}$ of the Weyl group,  where $\tau$ changes the sign of $x_6$, and $S_6$ acts by permuting the first six coordinates. Then the root $x_1 + \dotsb + x_5 - x_6 + x_7 - x_8$ is invariant under $H$ and contained in the connected hyperplane $E(M) \cap \langle x_1, x_2, x_3, x_4, x_5 - x_6, x_7 - x_8 \rangle$, whose orbit under $H$ has size 15.  This shows our claim.

For $M$ of type $E_6$, we argue in an analogous way: every element has $15$ rank one neighbors, 
and the maximal cardinality of a connected hyperplane is the cardinality of the $D_5$ matroid, i.e. $20$. We need to show that every element has at least 6 corank one neighbors. As above, we conclude with the observation that the element $x_1 + \dotsb + x_5 - x_6 - x_7 + x_8$ is contained in the connected hyperplane $E(M) \cap \langle  x_1, x_2, x_3, x_4 + x_5, -x_6 - x_7 + x_8 \rangle$, 
whose orbit under $S_5$ has size 10. 

The rank $3$ matroid associated to the Coxeter arrangement of type $H_3$ does not admit a Cremona basis (see e.g. the list of flats in \cite{H3}), hence in this case our claim follows from \cite[Theorem 9.5]{sw}.

The matroid $M$ associated to a Coxeter arrangement of type $H_4$ can be realized by the matrix $H$ over $\mathbb{Q}[\frac{1}{2}(1 + \sqrt{5})]$ given in \cite[Definition 2.1]{H4}. It is simple of rank $4$ and consists of $60$ elements. 
By Lemma \ref{flats} every rank $3$ flat in $M$ corresponds to a Coxeter arrangement given by a subgraph of the Coxeter graph of $H_4$ with three vertices. Hence every rank $3$ flat in $M$ has at most $15$ elements, which is the cardinality of the matroid associated to an $H_3$ type Coxeter arrangement.
By 
\cite[Proposition 3.5(4)]{H4}, every element in $E(M)$ forms a disconnected flat together with at least $15$ other elements contained in orthoframes, so that every rank one ray has at least $15$ rank one neighbours. Since it is contained in at least one connected hyperplane, the rank $1$ vertices of $S$ have strictly more neighbors than the corank $1$ vertices of $S$.

\end{proof}


\vspace{3ex}
\author{Stefan Rettenmayr,}
\address{Goethe University Frankfurt, Institut f\"ur Mathematik, Robert-Mayer-Strasse 6-8, 60325 Frankfurt, Germany.}
\email{rettenmayr@math.uni-frankfurt.de} 

\vspace{1ex}
\noindent \author{Annette Werner,}
\address{Goethe University Frankfurt, Institut f\"ur Mathematik, Robert-Mayer-Strasse 6-8, 60325 Frankfurt, Germany.}
\email{werner@math.uni-frankfurt.de}

\end{document}